\documentclass[12pt]{amsart}
\usepackage{amsmath, amssymb, latexsym, amsthm, mathrsfs, color}
\usepackage[all]{xy}
\long\def\forget#1\forgotten{}

      \newenvironment{changemargin}[2]{\begin{list}{}{
         \setlength{\topsep}{0pt}\setlength{\leftmargin}{0pt}
         \setlength{\rightmargin}{0pt}
         \setlength{\listparindent}{\parindent}
         \setlength{\itemindent}{\parindent}
         \setlength{\parsep}{0pt plus 1pt}
         \addtolength{\leftmargin}{#1}\addtolength{\rightmargin}{#2}
         }\item }{\end{list}}

\newcommand{\nc}{\newcommand}
\nc{\card}[1]{{\bigl|#1\bigr|}}
\nc{\set}[2]{\{#1\,:\,#2\}}
\nc{\varset}[2]{\left\{#1\,:\,#2\right\}}
\nc{\bt}{\op{bt}}
\nc{\bts}{\op{bts}}
\nc{\UN}{\op{U}}
\nc{\aL}{\op{aL}}
\nc{\fbt}{\op{fbt}}
\nc{\N}{\mathbb{N}}
\nc{\cN}{\mathcal{N}}
\nc{\cU}{\mathcal{U}}
\nc{\cV}{\mathcal{V}}
\nc{\Nb}{\cN}
\nc{\cNb}{\cN_\theta}
\nc{\ct}{\chi_\theta}
\nc{\cF}{\mathcal{F}}
\nc{\thd}[1]{\theta(#1)}
\nc{\thh}[1]{\cl{#1}^\theta}
\nc{\my}[1]{{\color{red}[#1]}}
\nc{\Pa}[8]{\bibitem{#1} {#2}, \emph{#3}, {#4} \textbf{#5} ({#6}), {#7}--{#8}.}
\nc{\tPa}[5]{\bibitem{#1} {#2}, \emph{#3}, {#4}, to appear.}
\nc{\sPa}[4]{\bibitem{#1} {#2}, \emph{#3}, {#4}, submitted.}
\nc{\Bc}[9]{\bibitem{#1} {#2}, \emph{#3}, in: \textbf{#4} (#5), #6 #7, #8--#9.}
\nc{\mx}[1]{\begin{matrix}#1\end{matrix}}
\nc{\arx}[1]{\texttt{http://arxiv.org/math/#1}}
\nc{\bq}{\begin{quote}}
\nc{\eq}{\end{quote}}
\nc{\cl}[1]{\overline{#1}}
\nc{\alephes}{{\aleph_0}}
\nc{\op}{\operatorname}
\nc{\Union}{\bigcup}
\long\def\forget#1\forgotten{}
\nc{\oo}{\infty}
\nc{\x}{\times}
\nc{\sub}{\subseteq}
\nc{\spst}{\supseteq}
\nc{\sm}{\setminus}

\nc{\la}{\langle}
\nc{\ra}{\rangle}
\newtheorem{thm}{Theorem}[section]\nc{\bthm}{\begin{thm}} \nc{\ethm}{\end{thm}}
\newtheorem{prop}[thm]{Proposition}\nc{\bprp}{\begin{prop}} \nc{\eprp}{\end{prop}}
\newtheorem{fact}[thm]{Fact}\nc{\bfct}{\begin{fact}} \nc{\efct}{\end{fact}}
\newtheorem{prob}[thm]{Problem}\nc{\bprb}{\begin{prob}} \nc{\eprb}{\end{prob}}
\newtheorem{lem}[thm]{Lemma}\nc{\blem}{\begin{lem}} \nc{\elem}{\end{lem}}
\newtheorem{claim}[thm]{Claim}\nc{\bclm}{\begin{claim}} \nc{\eclm}{\end{claim}}
\newtheorem{cor}[thm]{Corollary}\nc{\bcor}{\begin{cor}} \nc{\ecor}{\end{cor}}
\newtheorem{conj}[thm]{Conjecture}\nc{\bcnj}{\begin{conj}} \nc{\ecnj}{\end{conj}}
\theoremstyle{definition}
\newtheorem{defn}[thm]{Definition}\nc{\bdfn}{\begin{defn}} \nc{\edfn}{\end{defn}}
\theoremstyle{remark}
\newtheorem{rem}[thm]{Remark}\nc{\brem}{\begin{rem}} \nc{\erem}{\end{rem}}
\newtheorem{cnv}[thm]{Convention}\nc{\bcnv}{\begin{cnv}} \nc{\ecnv}{\end{cnv}}
\newtheorem{exam}[thm]{Example}\nc{\bexm}{\begin{exam}} \nc{\eexm}{\end{exam}}
\nc{\bpf}{\begin{proof}} \nc{\epf}{\end{proof}}
\nc{\be}{\begin{enumerate}}\nc{\ee}{\end{enumerate}}
\nc{\bi}{\begin{itemize}}\nc{\itm}{\item}\nc{\ei}{\end{itemize}}
\nc{\Subsection}[1]{\goodbreak\subsection*{#1}}

\forget
\forgotten

\forget
\forgotten

\newcommand{\ed}{

\subsection*{Final comment}
Replacing, everywhere relevant, \emph{closed neighborhoods} by \emph{neighborhoods},
one obtains the notions of \emph{finitely-Hausdorff} spaces, and the corresponding results hold true.
This line of investigation was initiated by Bonanzinga in \cite{Bo}.
The results presented here generalize some of her results.

\subsection*{Acknowledgments}
This research was partially supported by CNR (GNSAGA) and MIUR, Italy, through ``Fondi 40\%''.
A part of this work was carried out during a visit of the fourth named author at the
University of Messina.
This author thanks his hosts for their warm hospitality and stimulating atmosphere.

\end{document}}

\title[The $\theta$-closed hull]{On the cardinality of the\\ $\theta$-closed hull of sets}
\author[Cammaroto]{Filippo Cammaroto}
\author[Catalioto]{Andrei Catalioto}
\author[Pansera]{Bruno Antonio Pansera}
\address[Filippo Cammaroto, Andrei Catalioto, Bruno Antonio Pansera]{Dipartimento di Matematica,
Universit\'a di Messina, viale Ferdinando Stagno D'Al\-contres 31, S. Agata 98166, Messina, Italy}
\email{[camfil,acatalioto,bpansera]@unime.it}

\author[Tsaban]{Boaz Tsaban}
\address[Boaz Tsaban]{Department of Mathematics, Bar-Ilan University, Ram\-at Gan 52900, Israel}
\email{tsaban@math.biu.ac.il}
\urladdr{http://www.cs.biu.ac.il/\~{}tsaban}

\keywords{Urysohn space, $n$-Urysohn space, finitely-Urysohn space, Urysohn number,
H-closed space, H-set,
$\theta$-closure, $\theta$-closed hull, $\theta$-tightness, $\theta$-bitightness,
finite $\theta$-bitightness, $\theta$-bitightness small number,
$\theta$-character, character, cardinal inequalities.}

\subjclass{%
Primary: 54A25, 54D10; Secondary: 54A20, 54D25.
}

\begin{document}

\begin{abstract}
The \emph{$\theta$-closed hull} of a set $A$ in a topological space is the smallest
set $C$ containing $A$ such that, whenever all \emph{closed} neighborhood of a point intersect
$C$, this point is in $C$.

We define a new topological cardinal invariant function, the
\emph{$\theta$-bitighness small number} of a space $X$, $\bts_\theta(X)$, and prove that
in every topological space $X$, the cardinality of the $\theta$-closed hull
of each set $A$ is at most $\card{A}^{\bts_\theta(X)}$.
Using this result, we synthesize all earlier results on bounds on
the cardinality of $\theta$-closed hulls.
We provide applications to P-spaces and to the almost-Lindel\"of number.
\end{abstract}

\maketitle

\section{Introduction}

An \emph{Urysohn} (or \emph{T$_{2\frac{1}{2}}$}) \emph{space},
is a space in which distinct points are separated by closed neighborhoods.
Thus, Urysohn spaces are in between Hausdorff and regular spaces.
The spaces considered here generalize Urysohn spaces.

Let $X$ be a topological space.
A point $x\in X$ is in the \emph{$\theta$-derivative} $\thd{A}$ of a set $A\sub X$
if each closed neighborhood of $x$ intersects $A$ (cf.\ Veli\v{c}ko \cite{Velicko69}).\footnote{The
use of the letter $\theta$ for these concepts was proposed by Alexandroff, in recognition of
Fedorchuk's results on the involved concepts
($\theta$ is the first letter in a Greek transcription of ``Fedorchuk''.)
For additional details on the history of and motivation for the concepts treated in this paper, see \cite{CGNP}.}
For regular spaces, $\thd{A}=\cl{A}$,
but in general the operator $\theta$ is not idempotent for Urysohn spaces.\footnote{In
earlier works, the $\theta$-derivative $\thd{A}$ is also denoted $\op{cl}_\theta(A)$ and called
\emph{$\theta$-closure}. Since the operator $\theta$ is not idempotent, we decided not to use the term
\emph{closure} here.}
The \emph{$\theta$-closed hull} $\thh{A}$ of $A$ (cf.\ \cite{Be-Ca}) is
the smallest set $C\sub X$ such that $A\sub C=\theta(C)$.\footnote{In
earlier works, the $\theta$-closed hull of $A$ is also denoted $[A]_\theta$.}

As there are first countable Urysohn spaces $X$ and sets $A\sub X$ such that, e.g.,
$\card{\cl{A}}=\alephes<2^\alephes=\thd{A}$ \cite{Be-Ca}, a major goal concerning
the mentioned concepts is
that of providing upper bounds on the cardinalities of $\theta$-closed hulls of sets,
in terms of cardinal functions of the ambient space $X$
(e.g., Bella--Cammaroto \cite{Be-Ca}, Cammaroto--Ko\v{c}inac \cite{Ca-Ko1, Ca-Ko2}, Bella \cite{Be},
Alas--Ko\v{c}inac \cite{Al-Ko}, Bonanzinga--Cammaroto--Matveev \cite{Bo-Ca-Ma}, Bonanzinga--Pansera
\cite{Bo-Pa}, and McNeill \cite{Mc}).
We identify several concepts and topological cardinal functions,
which lead to generalizations of results from the mentioned papers.

\medskip

Throughout this paper, $X$ is a topological space and $A$ is an arbitrary subset of $X$.
\medskip

Recall that for $x\in X$, $\chi(X,x)$ is the minimal cardinality of a local base at $x$,
and the \emph{character} $\chi(X)$ of $X$ is the maximum of $\alephes$ and $\sup_{x\in X}\chi(X,x)$.
In 1988, Bella and Cammaroto proved that, for Urysohn spaces $X$, $\card{\thh{A}}\le\card{A}^{\chi(X)}$ \cite{Be-Ca}.

For $x\in X$, let $\ct(X,x)$ be the minimal cardinality of a family of \emph{closed} neighborhoods of $x$
such that each closed neighborhood of $x$ contains one from this family.
The \emph{$\theta$-character} $\ct(X)$ of $X$ is the maximum of $\alephes$ and $\sup_{x\in X}\ct(X,x)$.
Thus, $\ct(X)\le\chi(X)$.
In \cite{Al-Ko}, Alas and Ko\v{c}inac define this topological cardinal invariant,
show that the inequality may be proper, and modify the Bella--Cammaroto argument to
show that, for Urysohn spaces $X$, $\card{\thh{A}}\le\card{A}^{\ct(X)}$.

In 1993, Cammaroto and Ko\v{c}inac defined the \emph{$\theta$-bitighness} of an Urysohn space $X$,
$\bt_\theta(X)$, to be the minimal cardinal $\kappa$ such that,
for each non-$\theta$-closed $A\sub X$,  there are $x\in\thd{A}\sm A$
and sets $A_\alpha \in [A]^{\leq \kappa}$, $\alpha<\kappa$, such that $\bigcap_{\alpha<\kappa}\thd{A_\alpha}=\{x\}$
\cite{Ca-Ko1}.
For Urysohn spaces $X$, Cammaroto and Ko\v{c}inac proved that $\bt_\theta(X)\le\chi(X)$.
Moreover, their proof shows that $\bt_\theta(X)\le\ct(X)$.
They supplied examples where the inequality is strict, and proved that $\card{\thh{A}}\le\card{A}^{\bt_\theta(X)}$,
thus refining the Bella--Cammaroto Theorem.

In their recent work \cite{Bo-Ca-Ma}, Bonanzinga, Cammaroto and Matveev defined the \emph{Urysohn number}
$\UN(X)$ to be the minimal cardinal $\kappa$ such that, for each set $\set{x_\alpha}{\alpha<\kappa}\sub X$,
there are closed neighborhoods $U_\alpha$ of $x_\alpha$, $\alpha<\kappa$,
such that $\bigcap_{\alpha<\kappa}U_\alpha=\emptyset$.
Thus, $X$ is Urysohn if and only if $\UN(X)=2$.
They note that, for Hausdorff spaces, $\UN(X)\le |X|$, and prove that for each cardinal $\kappa\ge 2$,
there is a Hausdorff space with $\UN(X)=\kappa$ \cite{Bo-Ca-Ma}.

\bdfn
$X$ is \emph{finitely-Urysohn} if $\UN(X)$ is finite.
\edfn

Bonanzinga, Cammaroto and Matveev generalized Bella and Cammaroto's result from Urysohn
to finitely-Urysohn spaces \cite{Bo-Ca-Ma}.
Later, Bonanzinga and Pansera improved this and the Alas--Ko\v{c}inac result:
For finitely-Urysohn spaces, $\card{\thh{A}}\le\card{A}^{\ct(X)}$ \cite{Bo-Pa}.

A technical problem in synthesizing the Bonanzinga--Cammaroto--Matveev Theorem
and the the Cammaroto--Ko\v{c}inac Theorem is that $\bt_\theta(X)$ need not be
defined for finitely-Urysohn spaces.

We define a new topological cardinal invariant function, the
\emph{$\theta$-bitighness small number} of a space $X$, denoted $\bts_\theta(X)$, and prove the following
assertions:
\be
\item $\bts_\theta(X)$ is defined for all topological spaces $X$ (Definition \ref{btsdef}).
\item Whenever $\bt_\theta(X)$ is defined, $\bts_\theta(X)\le \bt_\theta(X)$ (Corollary \ref{immediate} 
and Definition \ref{btsdef}).
\item For all finitely-Urysohn spaces, $\bts_\theta(X)\le\ct(X)$ (Theorem \ref{fbtLTct} and Definition \ref{btsdef}).
\item In every topological space $X$, $\card{\thh{A}}\le\card{A}^{\bts_\theta(X)}$ (Theorem \ref{btsthm}).
\ee
This generalizes all of the above-mentioned results.
The situation is summarized in the following diagram.
\begin{changemargin}{-5cm}{-5cm}
\par\bigskip
\begin{center}
$\xymatrix{
\forall\,\mbox{finitely-Urysohn }X,\ \card{\thh{A}}\le\card{A}^{\chi(X)} \mbox{ \cite{Bo-Ca-Ma}}\ar[r]
& \forall\,\mbox{Urysohn }X,\ \card{\thh{A}}\le\card{A}^{\chi(X)} \mbox{ \cite{Be-Ca}}\\
\forall\,\mbox{finitely-Urysohn }X,\ \card{\thh{A}}\le\card{A}^{\ct(X)} \mbox{ \cite{Bo-Pa}}\ar[r]\ar[u]
& \forall\,\mbox{Urysohn }X,\ \card{\thh{A}}\le\card{A}^{\ct(X)} \mbox{ \cite{Al-Ko}}\ar[u]\\
\forall\,X,\ \card{\thh{A}}\le\card{A}^{\bts_\theta(X)}\ar[r]\ar[u]
& \forall\,\mbox{Urysohn }X,\ \card{\thh{A}}\le\card{A}^{\bt_\theta(X)} \mbox{ \cite{Ca-Ko1}}\ar[u]
}$
\end{center}
\par\bigskip
\end{changemargin}
We actually establish finer theorems than the ones mentioned above, as explained in the
following sections.

We also provide a partial solution to a problem of Bonanzinga--Cammaroto--Matveev \cite{Bo-Ca-Ma} and
Bonanzinga--Pansera \cite{Bo-Pa}.

\section{Finite bitightness and the bitightness small number}

\bdfn
The \emph{finite $\theta$-bitighness} of a space $X$,
$\fbt_\theta(X)$, is the minimal cardinal $\kappa$ such that,
for each non-$\theta$-closed $A\sub X$,  there are sets $A_\alpha \in [A]^{\leq \kappa}$, $\alpha<\kappa$,
such that $\bigcap_{\alpha<\kappa}\thd{A_\alpha}\sm A$ is finite and nonempty.
\edfn

\bcor\label{immediate}
$\fbt_\theta(X)$ is defined for all finitely-Urysohn spaces.
When $\bt_\theta(X)$ is defined, so is $\fbt_\theta(X)$, and $\fbt_\theta(X)\le\bt_\theta(X)$.\qed
\ecor

The following easy fact will be used in several occasions.

\blem\label{yey}
If $x\in\thd{A}$, then for each closed neighborhood $V$ of $x$, $x\in\thd{A\cap V}$.\qed
\elem

For Urysohn spaces, $\fbt_\theta(X)$ is very closely related to $\bt_\theta(X)$.

\bprp
Let $X$ be an Urysohn space, and $\kappa=\fbt_\theta(X)$.
For each non-$\theta$-closed $A\sub X$, there are
$x\notin A$ and $A_\alpha \in [A]^{\leq \kappa}$, $\alpha<\kappa$,
such that $\bigcap_{\alpha<\kappa}\thd{A_\alpha}\sm A=\{x\}$.
\eprp
\bpf
Pick sets $A_\alpha \in [A]^{\leq \kappa}$, $\alpha<\kappa$,
such that $\bigcap_{\alpha<\kappa}\thd{A_\alpha}\sm A$ is finite, say equal to $\{x_1,\dots,x_k\}$.

Since $X$ is Urysohn, there are closed neighborhoods $V_i$ of $x_i$, $i\le k$, such that
$V_1\cap (V_2\cup\cdots\cup V_k)=\emptyset$. Indeed, for each $i=2,\dots,k$ pick disjoint closed
neighborhoods $U_i$ and $V_i$ of $x_1,x_i$, respectively, and set $V_1=U_2\cap\cdots\cap U_k$.

For each $\alpha<\kappa$, $x_1\in\thd{A_\alpha\cap V_1}$.
Then $A_\alpha\cap V_1\in [A]^{\leq \kappa}$ for each $\alpha$, and
$$\bigcap_{\alpha<\kappa}\thd{A_\alpha\cap V_1}\sm A=\{x_1\}.\qedhere$$
\epf

\blem\label{finlemma}
Let $X$ be a finitely-Urysohn space.
For all $B,D \subseteq X$ with $B\sub\thd{D}$ and $|B| \ge \UN(X)$, there are $1\le m\le k<\UN(X)$
and $b_1,\dots, b_k \in B$ such that
\begin{equation}\label{int}
B \cap \bigcap_{V \in \cNb(b_1) \wedge \cdots \wedge \cNb(b_k)} \thd{D\cap V} = \{b_1,\dots,b_m\}.
\end{equation}
\elem
\bpf
For $k=\UN(X)$, any intersection as in \eqref{int} is empty.
For $k=1$, any such intersection is nonempty (since, by Lemma \ref{yey}, it contains $b_1$).
Thus, let $k$ be maximal such that there are $b_1,\dots,b_k\in B$ for which
the intersection in \eqref{int} is nonempty.
$1\le k<\UN(X)$.
We claim that
$$B\cap\bigcap_{V\in\cNb(b_1)\wedge\cdots\wedge\cNb(b_k)}\thd{D\cap V}\sub\{b_1,\dots,b_k\}.$$
Assume, towards a contradiction, that there is
$$x\in B\cap\bigcap_{V\in\cNb(b_1)\wedge\cdots\wedge\cNb(b_k)}\thd{D\cap V}\sm\{b_1,\dots,b_k\}.$$
By Lemma \ref{yey}, for each $V\in\cNb(b_1)\wedge\cdots\wedge\cNb(b_k)$ and each $W\in\cNb(x)$,
$x\in\thd{D\cap V\cap W}$. Thus,
$$x\in B\cap\bigcap_{V\in\cNb(b_1)\wedge\cdots\wedge\cNb(b_k)\wedge\cNb(x)}\thd{D\cap V},$$
and in particular this set is nonempty. This contradicts the maximality of $k$.

Thus, the intersection is nonempty, and by reordering $b_1,\dots,b_k$, we may assume
that the intersection is $\{b_1,\dots,b_k\}$ for some $m$ with $1\le m\le k$.
\epf

\bthm\label{fbtLTct}
For each finitely-Urysohn space $X$, $\fbt_\theta(X)\le\ct(X)$.
\ethm
\bpf
For families of sets $\cF_1, \cF_2, \dots, \cF_n \sub P(X)$, define
$$\cF_1 \wedge \cF_2 \wedge \cdots \wedge \cF_n :=
\set{V_1 \cap V_2 \cap \cdots \cap V_n}{V_1\in\cF_1,\dots,V_n\in\cF_n}.$$
For $x\in X$, let $\cNb(x)$ be the family of closed neighborhoods of $x$.

Let $\kappa=\ct(X)$. Let $A\sub X$ be non-$\theta$-closed.
Assume that $\thd{A}\sm A$ is finite. Fix $b\in \thd{A}\sm A$.
Fix a base $\set{V_\alpha}{\alpha<\kappa}$ for $\cNb(b)$.
For each $\alpha<\kappa$, let $a_\alpha\in A\cap V_\alpha$.
Let $D=\set{a_\alpha}{\alpha<\kappa}$, and set $A_\alpha=D$ for all $\alpha<\kappa$.
Then
$$b\in \thd{D}\sm A=\bigcap_{\alpha<\kappa} \thd{A_\alpha}\sm A\sub\thd{A}\sm A,$$
so that $\bigcap_{\alpha<\kappa} \thd{A_\alpha}\sm A$ is finite and nonempty,
and the requirement in the definition of $\fbt_\theta(X)\le\kappa$ is fulfilled.

Thus, assume that the set $B=\thd{A}\sm A$ is infinite.
Apply Lemma \ref{finlemma} to the sets $B$ and $D=A$, to obtain
$1\le m\le k<\UN(X)$ and $b_1,\dots,b_k\in B$ such that Equation \eqref{int} holds.
For each $i\le k$, fix a basis $\cF_i$ for $\cNb(b_i)$ with $|\cF_i|\le\kappa$.
Enumerate
$$\cF_1\wedge\cdots\wedge\cF_k=\set{V_\alpha}{\alpha<\kappa}.$$
By Equation \eqref{int},
$$B \cap \bigcap_{\alpha<\kappa} \thd{A\cap V_\alpha} = \{b_1,\dots,b_m\}.$$
In particular, for each $\alpha<\kappa$ there is $a_\alpha\in A\cap V_\alpha$.
Take $$C=\set{a_\alpha}{\alpha<\kappa}\in[A]^{\le\kappa}.$$
Fix $i\le m$ and $\alpha<\kappa$. Let $V\in\cNb(b_i)$.
Then $V_\alpha\cap V\in\cNb(b_1)\wedge\cdots\cNb(b_m)$, and thus there is
$\beta<\kappa$ such that $V_\beta\sub V_\alpha\cap V$.
Then $a_\beta\in C\cap V_\alpha\cap V$, and in particular $C\cap V_\alpha\cap V$ is nonempty.
This shows that $b_i\in\thd{C\cap V_\alpha}$.

Thus,
\begin{eqnarray*}
b_1,\dots,b_m & \in & \bigcap_{\alpha<\kappa}\thd{C\cap V_\alpha}\sm A\sub
\bigcap_{\alpha<\kappa}\thd{A\cap V_\alpha}\sm A\sub\\
& \sub & (\thd{A}\sm A)\cap \bigcap_{\alpha<\kappa}\thd{A\cap V_\alpha}=
\{b_1,\dots,b_m\},
\end{eqnarray*}
and therefore
$$\bigcap_{\alpha<\kappa}\thd{C\cap V_\alpha}\sm A=\{b_1,\dots,b_m\},$$
as required in the definition of $\fbt_\theta(X)\le\kappa$.
\epf


\bdfn\label{btsdef}
The \emph{$\theta$-bitightness small number} of $X$, $\bts_\theta(X)$, is
the minimal cardinal $\kappa$ such that, for each non-$\theta$-closed $A\sub X$ that is not a singleton,\footnote{In the Hausdorff context,
singletons are $\theta$-closed, and thus the restriction to non-singletons may be removed.}
there are $A_\alpha\in [A]^{\le\kappa}$, $\alpha<\kappa$, such that
$$\bigcap_{\alpha<\kappa}\thd{A_\alpha}\sm A\neq\emptyset\mbox{ and }
\card{\bigcap_{\alpha<\kappa}\thd{A_\alpha}}\le \card{A}^\kappa.$$
\edfn

$\bts_\theta(X)$ is defined for all spaces $X$, and is obviously $\le\fbt_\theta(X)$ whenever the latter is defined.

\bthm\label{btsthm}
Let $X$ be a topological space. For each $A\sub X$, 
$$\card{\thh{A}}\le\card{A}^{\bts_\theta(X)}.$$
\ethm
\bpf
Let $\kappa=\bts_\theta(X)$, $\lambda=\card{A}$.
We define sets $C_\alpha\sub X$, with $\card{C_\alpha}\le\lambda^\kappa$, $\alpha\le\kappa^+$, by induction on $\alpha$.

$C_0:=A$.

Given $C_\alpha$,
$$C_{\alpha+1}:=\Union\varset{\bigcap_{\beta<\kappa}\thd{B_\beta}}{\set{B_\beta}{\beta<\kappa}\sub [C_\alpha]^{\le\kappa},\
\card{\bigcap_{\beta<\kappa}\thd{B_\beta}}\le\lambda^\kappa}.$$
Then $C_\alpha\sub C_{\alpha+1}$. As $\card{C_\alpha}\le \lambda^\kappa$,
$\card{C_{\alpha+1}}\le ((\lambda^\kappa)^\kappa)^\kappa\cdot(\lambda^\kappa)^\kappa=\lambda^\kappa$.

For a limit ordinal $\alpha$, $C_\alpha:=\Union_{\beta<\alpha}C_\beta$.
Then $\card{C_\alpha}\le |\alpha|\cdot\lambda^\kappa\le\kappa^+\cdot\lambda^\kappa=\lambda^\kappa$.

End of the construction.

\medskip

Let $C=C_{\kappa^+}$. Then $\card{C}\le\lambda^\kappa$, $A=C_0\sub C$, and $C$ is $\theta$-closed.
Indeed, assume otherwise and let $B_\alpha\in [C]^{\le\kappa}$, $\alpha<\kappa$, be such that
$\bigcap_{\alpha<\kappa}\thd{B_\alpha}\sm C\neq\emptyset$ and
$\card{\bigcap_{\alpha<\kappa}\thd{B_\alpha}}\le \card{C}^\kappa$.
Then $\card{\bigcap_{\alpha<\kappa}\thd{B_\alpha}}\le (\lambda^\kappa)^\kappa=\lambda^\kappa$.
As $\kappa^+$ is regular, for each $\alpha<\kappa$ there is $\beta_\alpha<\kappa^+$ such that $B_\alpha\sub C_{\beta_\alpha}$.
Again as $\kappa^+$ is regular, $\beta:=\sup_{\alpha<\kappa}\beta_\alpha<\kappa$.
Then $B_\alpha\in [C_\beta]^{\le\kappa}$ for all $\alpha<\kappa$, and thus $\bigcap_{\alpha<\kappa}\thd{B_\alpha}\sub C_{\beta+1}\sub C$.
A contradiction.
\epf

\brem
Immediately after Proposition 7 of \cite{Be}, Bella points out that there are Hausdoff spaces $X$
where the inequality $\card{\thh{A}}\le\card{A}^{\chi(X)}$ fails for some of their subsets.
In particular, by Theorem \ref{btsthm},
$\bts_\theta(X)$ may be larger than $\ct(X)$ may fail for general Hausdorff spaces $X$.
\erem

\section{The $\theta$-closed hull in P-spaces}

Bonanzinga--Cammaroto--Matveev \cite{Bo-Ca-Ma} and Bonanzinga--Pansera \cite{Bo-Pa} ask whether, in
all Hausdorff spaces $X$, $\card{\thh{A}}\le\card{A}^{\ct(X)}\cdot\UN(X)$. We give a partial answer.

\bdfn
The \emph{$\theta$-P-point number} of a space is the minimal cardinal $\kappa$ such that some $x\in X$ has
closed neighborhoods $V_\alpha$, $\alpha<\kappa$, with $\bigcap_{\alpha<\kappa}V_\alpha$ not a neighborhood of $x$.
\edfn

As the \emph{$\theta$-P-point number} of any space is at leat $\alephes$, the following
theorem generalizes the Bonanzinga--Pansera Theorem, and thus also the earlier three theorems
discussed in the introduction.

\bthm
Let $X$ be a topological space whose Urysohn number is smaller than its $\theta$-P-point number.
For each $A\sub X$, 
$$\card{\thh{A}}\le\card{A}^{\ct(X)}\cdot\UN(X).$$
\ethm
\bpf
Let $\kappa=\ct(X)$. For each $x\in\thd{A}$, let $\set{V^x_\alpha}{\alpha<\kappa}$ be a family of
closed neighborhoods of $x$ such that each closed neighborhood of $x$ contains one from this family.
For each $\alpha<\kappa$, fix $a_{x,\alpha}\in A\cap V^x_\alpha$. Let $A_x=\set{a_{x,\alpha}}{\alpha<\kappa}$.

Define a map
\begin{eqnarray*}
\Psi\colon  \thd{A} & \to & [[A]^{\le\kappa}]^{\le\kappa}\\
            x & \mapsto & \set{A_x\cap V^x_\alpha}{\alpha<\kappa}.
\end{eqnarray*}
Let $\nu=\UN(X)$.
Let $x_\alpha$, $\alpha<\nu$, be distinct elements of $\thd{A}$ which are all mapped
to the same element $\Psi(x)$.
For each $\alpha<\nu$, pick $\beta_\alpha<\kappa$ such that
$$\bigcap_{\alpha<\nu}V^{x_\alpha}_{\beta_\alpha}=\emptyset.$$
Let $\alpha<\nu$. As $\Psi(x_\alpha)=\Psi(x)$, there is
$\gamma_\alpha<\kappa$ such that $A_{x_\alpha}\cap V^{x_\alpha}_{\beta_\alpha}=A_x\cap V^x_{\gamma_\alpha}$.
As $\nu$ is smaller than the $\theta$-P-point number of $X$, $V:=\bigcap_{\alpha<\nu}V^x_{\gamma_\alpha}$
is a closed neighborhood of $x$. Fix $\delta<\kappa$ such that $V^x_\delta\sub V$.
Then
\begin{eqnarray*}
a_{x,\delta} & \in & A_x\cap V^x_\delta\sub A_x\cap V=A_x\cap \bigcap_{\alpha<\nu}V^x_{\gamma_\alpha}=
\bigcap_{\alpha<\nu}A_x\cap V^x_{\gamma_\alpha}=\\
& = & \bigcap_{\alpha<\nu}A_{x_\alpha}\cap V^{x_\alpha}_{\beta_\alpha}\sub \bigcap_{\alpha<\nu} V^x_{\beta_\alpha}=\emptyset;
\end{eqnarray*}
a contradiction.

Thus, $\Psi$ is $<\nu$ to $1$, and therefore the cardinality of $\thh{A}$ is at most
$$\card{[[A]^{\le\kappa}]^{\le\kappa}}\cdot\nu= \card{A}^\kappa \cdot\nu.$$

By induction on $\alpha\le\kappa^+$, define $A_0:=A$, $A_{\alpha+1}:=\thd{A_\alpha}$, and
$A_\alpha=\Union_{\beta<\alpha}A_\beta$ for limit ordinals $\alpha$.
Then, by induction, $\card{A_\alpha}\le \card{A}^\kappa \cdot\nu$ for all $\alpha$.
As $\ct(X)=\kappa$, $A_{\kappa^+}=\thh{A}$ \cite{Bo-Pa}.
\epf

Recall that $X$ is a \emph{P-space} if each countable intersection of neighborhoods is a neighborhood.
Thus, the $\theta$-P-point number of a P-space is $\ge\aleph_1$.

\bcor
Let $X$ be a P-space with $\UN(X)=\alephes$.
For each $A\sub X$, $\card{\thh{A}}\le\card{A}^{\ct(X)}$.\qed
\ecor

\section{The almost-Lindel\"of number}

\bdfn[\cite{Be-Ca}]
The \emph{almost-Lindel\"of number} $\aL(A,X)$ of a set $A\sub X$ is the minimal cardinal $\kappa$ such that,
for each open cover $\cU$ of $A$, there is $\cV\in[\cU]^{\le\kappa}$ such that $A\sub\Union_{U\in\cU}\cl{U}$.
\edfn

\bthm\label{main2}
Let $X$ be a Hausdorff topological space.
For each $A\sub X$,
$$\card{A}\le 2^{\bts_\theta(X)\cdot\ct(X)\cdot\aL(A,X)}.$$
\ethm
\bpf
Let $\kappa=\bts_\theta(X)\cdot\ct(X)\cdot\aL(A,X)$.
For each $x\in X$, let $\cF_x$ be a family of closed neighborhoods of $x$ such that $|\cF_x|\le\kappa$,
and each closed neighborhood of $x$ contains one from $\cF_x$.

Fix $a\in A$. We define, by induction on $\alpha\le\kappa^+$, sets $A_\alpha\sub X$ such that $\card{A_\alpha}\le2^\kappa$.

$A_0:=\{a\}$.

Step $\alpha>0$: Let $B=\Union_{\beta<\alpha}A_\beta$. By the induction hypothesis, $\card{B}\le 2^\kappa$.
Thus, $\card{\Union_{x\in B\cap A}\cF_x}\le 2^\kappa$ as well, and therefore
$\card{[\Union_{x\in B\cap A}\cF_x]^{\le\kappa}}\le 2^\kappa$.
For each $\cV\in [\Union_{x\in B\cap A}\cF_x]^{\le\kappa}$, with $A\sm\Union\cV\neq\emptyset$,
pick a point from $A\sm\Union\cV$. Let $C$ be the set of these points.
Then $|B\cup C|\le 2^\kappa$.
Set $B_\alpha=\thh{B\cup C}$. As $\bts_\theta(X)\le\kappa$, we have by Theorem \ref{btsthm} that
$\card{B}\le (2^\kappa)^{\kappa}=2^\kappa$.
End of the construction.

Let $B=B_{\kappa^+}$.
It remains to show that $A\sub B$. Assume otherwise, and fix $a_0\in A\sm B$.
As $B$ is $\theta$-closed, for each $x\in A\sm B$ we can choose $V_x\in\cF_x$
such that $V_x\cap B=\emptyset$.
For $x\in A\cap B$, choose $V_x\in\cF_x$ such that $a_0\notin V_x$.
As $\set{{V_x}^\circ}{x\in X}$ is an open cover of $A$ and $A$ is an H$_\kappa$-set,
there is $K\in [X]^{\le\kappa}$ such that $A\sub \Union_{x\in K}V_x$.
As $V_x\cap  B=\emptyset$ for each $x\in A\sm B$,
$$B\cap A\sub \Union_{x\in K\cap B}V_x.$$
As $\kappa^+$ is regular, there is $\alpha<\kappa^+$ such that $K\cap B\sub B_\alpha$.
As $a_0\in A\sm \Union_{x\in K\cap B}V_x$, we have by the construction of $B_{\alpha+1}$
an element in $B_{\alpha+1}\cap A\sm \Union_{x\in K\cap B}V_x$, and therefore so in
$B\cap A\sm \Union_{x\in K\cap B}V_x$;
a contradiction.
\epf

The following corollary improves upon a result of Bonanzinga, Cammaroto and Matveev \cite{Bo-Ca-Ma},
asserting that for Hausdorff, finitely-Urysohn spaces $X$, $\card{X}\le 2^{\chi(X)\cdot\aL(X,X)}$.

\bcor
Let $X$ be a Hausdorff, finitely-Urysohn space. For each $A\sub X$, $\card{A}\le 2^{\ct(X)\cdot\aL(A,X)}$.
In particular, $\card{X}\le 2^{\ct(X)\cdot\aL(X,X)}$.
\ecor
\bpf
By Theorem \ref{fbtLTct}, $\bts_\theta(X)\le\fbt_\theta(X)\le\ct(X)$ for finitely-Urysohn spaces.
Thus, Theorem \ref{main2} applies.
\epf

\ed